\documentclass[a4paper,12pt]{amsart}
\usepackage{amssymb}
\usepackage{graphicx}
\usepackage{ifthen}
\nonstopmode \numberwithin{equation}{section}
\newtheorem{thm}{Theorem}
\newtheorem{cor}{Corollary}
\newtheorem{lem}{Lemma}

\newtheorem{conj}{Conjecture}

\theoremstyle{definition}
\newtheorem{defn}{Definition}[section]

\newtheorem{prob}[equation]{Problem}

\newenvironment{rem}{%
\bigskip
\noindent \textsl{{\sl Remark. }}}{\bigskip}
\newenvironment{rems}{%
\bigskip
\noindent \textsl{{\sl Remarks. }}}{\bigskip}

\newcounter {own}
\def\theown {\thesection       .\arabic{own}}

\newenvironment{pf}[1][]{%
 \vskip 3mm
 \noindent
 \ifthenelse{\equal{#1}{}}%
  {{\slshape Proof. }}%
  {{\slshape #1.} }%
 }%
{\qed\bigskip}

\newcounter{alphabet}

\newtheorem{theorem}{Theorem}

\newcommand{\IR}{{\mathbb R}}
\newcommand{\ID}{{\mathbb D}}
\newcommand{\IN}{{\mathbb N}}
\newcommand{\IC}{{\mathbb C}}

\newcommand{\M}{{\mathcal M}}

\newcommand{\IB}{{\mathcal B}}
\newcommand{\IH}{{\mathcal H}}

\def\be{\begin{equation}}
\def\ee{\end{equation}}

\newcommand{\bee}{\begin{enumerate}}
\newcommand{\eee}{\end{enumerate}}

\newcommand{\blem}{\begin{lem}}
\newcommand{\elem}{\end{lem}}
\newcommand{\bthm}{\begin{thm}}
\newcommand{\ethm}{\end{thm}}
\newcommand{\bcor}{\begin{cor}}
\newcommand{\ecor}{\end{cor}}
\newcommand{\beg}{\begin{examp}}
\newcommand{\eeg}{\end{examp}}
\newcommand{\begs}{\begin{examples}}
\newcommand{\eegs}{\end{examples}}
\newcommand{\bdefe}{\begin{defn}}
\newcommand{\edefe}{\end{defn}}
\newcommand{\bprob}{\begin{prob}}
\newcommand{\eprob}{\end{prob}}
\newcommand{\bei}{\begin{itemize}}
\newcommand{\eei}{\end{itemize}}

\newcommand{\bcon}{\begin{conj}}
\newcommand{\econ}{\end{conj}}
\newcommand{\bcons}{\begin{conjs}}
\newcommand{\econs}{\end{conjs}}
\newcommand{\bprop}{\begin{propo}}
\newcommand{\eprop}{\end{propo}}
\newcommand{\br}{\begin{rem}}
\newcommand{\er}{\end{rem}}
\newcommand{\brs}{\begin{rems}}
\newcommand{\ers}{\end{rems}}
\newcommand{\bo}{\begin{obser}}
\newcommand{\eo}{\end{obser}}
\newcommand{\bos}{\begin{obsers}}
\newcommand{\eos}{\end{obsers}}
\newcommand{\bpf}{\begin{pf}}
\newcommand{\epf}{\end{pf}}
\newcommand{\ba}{\begin{array}}
\newcommand{\ea}{\end{array}}
\newcommand{\beq}{\begin{eqnarray}}
\newcommand{\beqq}{\begin{eqnarray*}}
\newcommand{\eeq}{\end{eqnarray}}
\newcommand{\eeqq}{\end{eqnarray*}}

\newcommand{\ds}{\displaystyle}
\newcommand{\ov}{\overline}

\newcounter{minutes}
\divide\time by 60
\newcounter{hours}
\multiply\time by 60 \addtocounter{minutes}{-\time}
\begin{document}
\title{On the Bohr phenomenon for complex valued and vector valued functions}
\begin{center}
{\tiny \texttt{FILE:~\jobname .tex,
        printed: \number\year-\number\month-\number\day,
        \thehours.\ifnum\theminutes<10{0}\fi\theminutes}
}
\end{center}
\author{Bappaditya Bhowmik${}^{\mathbf{*}}$}
\address{Bappaditya Bhowmik, Department of Mathematics,
Indian Institute of Technology Kharagpur, Kharagpur - 721302, India.}
\email{bappaditya@maths.iitkgp.ac.in}
\author{Nilanjan Das}
\address{Nilanjan Das, Department of Mathematics,
Indian Institute of Technology Kharagpur, Kharagpur - 721302, India.}
\email{nilanjan@iitkgp.ac.in}

\keywords{Bohr inequality, compact group, Fourier transform, Hardy space.\newline
${}^{\mathbf{*}}$ Corresponding author}

\begin{abstract}
We explore the Bohr inequality involving the Fourier transforms of complex valued integrable and square integrable functions defined on a second countable compact topological group. We also investigate the connection of the Bohr phenomenon with a modulus of convexity of the space of bounded linear operators defined on a complex Hilbert space.
\end{abstract}
\thanks{The first author of this article would like to thank
SERB, DST, India (Ref.No.- MTR/2018/001176) for its financial support through MATRICS grant.}

\maketitle
\pagestyle{myheadings}
\markboth{B. Bhowmik, N. Das}{On the Bohr phenomenon for complex valued and vector valued functions}

\bigskip
\section{Introduction}
The famous Bohr's theorem (in improved form) reads as:

\begin{theorem}\label{TheoA}\cite{Bohr}
Let $f$ be a complex valued continuous function on the unit circle $\mathbb{T}$ with analytic Fourier series $\sum_{n=0}^\infty a_ne^{in\theta}$. If $f$ is bounded by $1$ on $\mathbb{T}$, then
\be\label{P6eq25}
\sum_{n=0}^\infty|a_ne^{in\theta}|r^n\leq 1
\ee
for any nonnegative real number $r$ satisfying $r\leq 1/3$. This constant $1/3$ is the best possible.
\end{theorem}
This theorem was originally proved for $r\leq 1/6$ by Harald Bohr, and it was further refined by Wiener, Riesz and Schur independently.
Moreover, this $1/3$ in Theorem \ref{TheoA} could be improved to $1/\sqrt{2}$ if $a_0=0$ (cf. \cite[Corollary 2.9]{Paul1}). Without loss of generality, one can consider $a_0\geq 0$ in Theorem \ref{TheoA}. With this assumption, it is shown in \cite[Theorem 2.1]{Paul1} that the boundedness hypothesis of $f$ can be weakened to $\mbox{Re}(f)\leq 1$. Using similar methods, Theorem \ref{TheoA} is further extended in the context of uniform algebras (see \cite{Paul2}).
Bohr obtained his result while attempting the absolute convergence problem for the Dirichlet series of the form $\sum a_nn^{-s}$, and it took several decades for the Bohr radius problem to emerge as an independent area of active research. In fact,
research on this topic picked up pace after Bohr's inequality was successfully applied to the characterization problem of Banach algebras satisfying the von Neumann
inequality (cf. \cite{Dix}).
Since then, Bohr's theorem continues to be studied in several different frameworks, viz. in multidimensional settings (see f.i. \cite{Aiz, Bay, Boas, Ga}), in certain abstract situations (cf. \cite{Ayt, Ham}), for ordinary and vector valued Dirichlet series (see for example \cite{Bala, Def2}), for a Faber-Green condenser (cf. \cite{La}), and for free holomorphic functions (see  \cite{Pop1}). The reader is also suggested to look at the references of the aforementioned articles to get a clearer picture of the recent developments on this subject.

We again concentrate on Theorem \ref{TheoA}. A closer examination of Wiener's proof of Theorem \ref{TheoA} in \cite[p. 4]{Bohr}, or the proof of Theorem 2.1 from \cite{Paul1} reveals that the key to establish Theorem \ref{TheoA} is the inequality $|a_n|\leq 2(1-a_0),\,n\geq 1$ (assuming $a_0\geq0$). Therefore, Bohr's theorem can be put in the following general form:
\begin{theorem}\label{TheoA1}
Under the hypotheses of Theorem \ref{TheoA} (with $a_0\geq 0$),
\be\label{P6eq26}
a_0+\sum_{n=1}^\infty|a_ne^{in\theta}z_n|\leq 1
\ee
for any sequence of complex numbers $\{z_n\}_{n=1}^\infty$ satisfying $\sum_{n=1}^\infty|z_n|\leq 1/2$.
\end{theorem}
From now onwards, we will use the term \emph{Bohr inequality} to denote inequalities of the type $(\ref{P6eq25})$ and $(\ref{P6eq26})$, as well as all of their variations. Presence of any such inequality in a result will be called the \emph{Bohr phenomenon} in this article.
It is now natural to ask if Theorem \ref{TheoA1} admits an extension for complex valued functions that do not necessarily have analytic Fourier series. We consider this question in a more general form.
In Theorem \ref{P6thm1} and Theorem \ref{P6thm2} of the present article, different versions of the Bohr inequality of type $(\ref{P6eq26})$ are proved for complex valued functions $f$ defined on a compact second countable topological group.

While the problem of finding various analogues of Bohr's theorem for the functions defined on a group has its own merits, Bohr inequalities for the vector valued functions defined on a group can also be useful
for characterizing a convexity property for certain complex Banach spaces. To illustrate this fact, we need some preparations. A complex Banach space $E$ with the norm $\|.\|_E$ is said to be \emph{$p$-uniformly $\IC$-convex} ($2\leq p<\infty$) if there exists a constant $\lambda(E)>0$ such that
\be\label{P6eq35}
(\|x\|_E^p+\lambda(E)\|y\|_E^p)^{1/p}\leq\max_\theta\|x+e^{i\theta}y\|_E
\ee
for all $x, y\in E$ (see f.i. \cite[Definition 1.9]{Bla2}). Suppose $H^\infty(\ID, E)$ is the space of bounded holomorphic functions $f$ from $\ID$ into $E$, and $\|f\|_{H^\infty(\ID, E)}=\sup_{|z|<1}\|f(z)\|_E$. It is known from \cite[Theorem 1.10]{Bla2} that a complex Banach space $E$ is $p$-uniformly $\IC$-convex ($2\leq p<\infty)$ if and only if there exists a constant $r_0(E)>0$ such that
$$
\left(\sum_{n=0}^\infty\|x_n\|_E^pr^{np}\right)^{\frac{1}{p}}\leq\|f\|_{H^\infty(\ID, E)}
$$
for $|z|=r\leq r_0(E)$, and for all $f(z)=\sum_{n=0}^\infty x_nz^n\in H^\infty(\ID, E)$. In the same spirit, Theorem \ref{P6thm3} of the present paper shows that $p$-uniform $\IC$-convexity of $\mathcal{B}(\IH)$, i.e. the space of bounded linear operators on a complex Hilbert space $\IH$, is equivalent to the existence of the Bohr phenomenon for a subclass of essentially bounded and strongly measurable operator valued functions defined on a compact abelian second countable group. It is easy to see that $\mathcal{B}(\IH)$ need not always be $p$-uniformly $\IC$-convex, take $\IH=\IC^2$ and
$$
x=
\begin{bmatrix}
1 & 0\\
0 & 0
\end{bmatrix}\,,
y=
\begin{bmatrix}
0 & 0\\
0 & 1
\end{bmatrix}.
$$
This paper is organized as follows. In section 2, we record
all the remaining prerequisites for the upcoming discussion. Section 3 contains the main results of this article and their proofs. The final section, i.e. section 4 includes some remarks in connection with the results obtained in this article.
\section{Notations and preliminaries}
We introduce a number of concepts and notations, which will be used frequently throughout this article.
Here and hereafter, $G$ denotes a second countable, compact topological group with a left and right invariant Haar measure $m$, normalized so that $m(G)=1$.
$G$ is not necessarily abelian unless mentioned specifically, and the topology of $G$ is always assumed to be Hausdorff.
A unitary representation of $G$ is a homomorphism $\pi$ from $G$ into the group $\mathcal{U}(\IH_\pi)$ of unitary operators on some nonzero Hilbert space $\IH_\pi$ that is continuous with respect to the strong operator topology, i.e. $x\mapsto\pi(x)u$ is continuous from $G$ to $\IH_\pi$ for any $u\in\IH_\pi$. $\IH_\pi$ is called the representation space of $\pi$, and its dimension is called the dimension or degree of $\pi$.
Let $\widehat{G}$ be the set of unitary equivalence classes $[\pi]$ of irreducible unitary representations of $G$. We clarify that two representations $\pi_1$ and $\pi_2$ are unitarily equivalent if there exists a unitary operator $U:\IH_{\pi_1}\to\IH_{\pi_2}$ satisfying
$\pi_2(x)=U\pi_1(x)U^{-1}$, $x\in G$. Moreover, a representation $\pi$ is said to be irreducible if there exists no nontrivial closed subspace $S$ of $\IH_\pi$ (i.e. $S\neq\{0\}$ and $\IH_\pi$) such that $\pi(x)S\subset S$ for all $x\in G$.
Since $G$ is both compact and second countable, i.e. has a countable basis as a topological space, from
\cite[Theorem (28.2), p. 61]{Hew} and \cite[Proposition (5.27), p. 138]{Foll} we conclude that $\widehat{G}$ is countable. We therefore write
$$
\widehat{G}=\{[\pi_n]: n\in\Lambda\cup\{0\}\}.
$$
Here $\pi_n$ is assumed to be the chosen representative for the class $[\pi_n]$, and $\Lambda$ is either whole of $\IN$ or a finite subset $\{1, 2, 3, \cdots, n_0\}$ of $\IN$.
Again, due to the compactness of $G$,
\cite[Theorem (5.2), p. 126]{Foll} asserts that
$$
d_{\pi_n}:=\dim(\IH_{\pi_n})<\infty.
$$
In our discussion, for any given complex Banach space $E$ with norm $\|.\|_E$, $L^p(G, E)$
is the space of strongly measurable functions $f$ from $G$ to $E$, satisfying $\|f\|_{L^p(G, E)}<\infty$. Here
\beqq \|f\|_{L^p(G, E)}&:=&\left(\int_G\|f(x)\|_E^p dx\right)^{1/p} \mbox {for}\, 1\leq p<\infty,\, \mbox {and}\\
 \|f\|_{L^\infty(G, E)}&:=& \ds\mbox{ess sup}_{x\in G}\|f(x)\|_E.
\eeqq
For the definition of strongly measurable functions and the related concept of Bochner integrability, we refer to \cite[p. 443]{Cue}. In particular, if $E=\IC$ then $\|.\|_E$ is the usual absolute value of complex numbers. Now, if $f\in L^1(G, \IC)$, the \emph{Fourier transform} of $f$ at $\pi_n$ is given by the operator (cf. \cite[p. 134]{Foll})
\be\label{P6eq6}
\widehat{f}(\pi_n):=\int_G f(x)\pi_n(x)^*\,dx
\ee
on $\IH_{\pi_n}$, where $\pi_n(x)^*$ is the adjoint of $\pi_n(x)$ for each $x\in G$. Considering an orthonormal basis $\{e_n\}_{n=1}^{d_{\pi_n}}$
for $\IH_{\pi_n}$, we can construct the matrix representation of $\pi_n(x)$. The $(i,j)$-th element of that representation is given by
$$
\pi_n(x)_{ij}=\langle\pi(x)e_j, e_i\rangle,
$$
where $\langle . , .\rangle$ is the inner product on $\IH_{\pi_n}$. Therefore, the $(i,j)$-th element of the matrix representation of $\widehat{f}(\pi_n)$ is
\be\label{P6eq8}
\widehat{f}(\pi_n)_{ij}=\int_G f(x)\ov{\pi_n(x)_{ji}}\,dx.
\ee
We will be working with these complex matrix forms of \lq$\widehat{f}(\pi_n)$\rq s and \lq$\pi_n$\rq s in this article.
Here we point out that, given any representation $\pi_n$ of $G$ on $\IH_{\pi_n}$ as above, there is another representation $\ov{\pi_n}$ (called the \emph{contragredient} of $\pi_n$) on the dual space of $\IH_{\pi_n}$. The $(i, j)$-th element of the matrix representation of $\ov{\pi_n}(x)$ is given by $\ov{\pi_n(x)_{ij}}$ (cf. \cite[p. 69]{Foll}).
Further, if we assume that $f\in L^2(G,\IC)$, then we have the following ``inversion formula":
$$
f(x)=\sum_{[\pi_n]\in\widehat{G}}d_{\pi_n}\mbox{tr}(\widehat{f}(\pi_n)\pi_n(x)),
$$
where the series at the right hand side converges in the $L^2$ norm.
We mention that for any complex matrix $A$, $\mbox{tr}(A)$ means the trace of the matrix $A$, i.e. sum of the diagonal entries of $A$, or equivalently, sum of the eigenvalues of $A$.
Also, we have
\be\label{P6eq9}
\|f\|_{L^2(G, \IC)}^2=\sum_{[\pi_n]\in\widehat{G}}d_{\pi_n}\mbox{tr}(\widehat{f}(\pi_n)^*\widehat{f}(\pi_n)).
\ee
For any $G$ we always have the trivial one dimensional irreducible representation which we call $\pi_0$, given by $\pi_0(x)=1$ for all $x\in G$. Starting from here, $\IH_{\pi_0}$ is taken to be $\IC$ in particular throughout this paper, and $d_{\pi_0}=1$. The reader might refer to \cite[Chapters 3, 5]{Foll} for an extensive treatment of the topics mentioned so far.

We will also be using a lot of matrix analysis, for which we need to be familiar with the following concepts.
For any $n\in\IN$, let $\M_n(\IC)$ be the space of $n$-dimensional matrices over $\IC$.
For any $i\in\IN$, $E_{ii}^{(n)}$ denotes the $d_{\pi_n}$ dimensional matrix with $(i,i)$-th entry $1$ and all other entries $0$.
For any two $A$, $B$ in $\M_n(\IC)$, the \emph{Hadamard product} $A\circ B$ is defined by
$$
(A\circ B)_{ij}=A_{ij}B_{ij},
$$
where $C_{ij}$ means the $(i,j)$-th entry of any $C\in\M_n(\IC)$.
Given any two $A, B\in\M_n(\IC)$, we will use the notation $A\boldsymbol{\cdot}B$ to denote both the usual matrix product $AB$ and the Hadamard product $A\circ B$.
Further, we fix another convention that for any $A, B, C\in\M_n(\IC)$, either  $A\boldsymbol{\cdot}B\boldsymbol{\cdot}C
=A\boldsymbol{\cdot}(B\boldsymbol{\cdot}C)
$,
or $A\boldsymbol{\cdot}B\boldsymbol{\cdot}C
=(A\boldsymbol{\cdot}B)\boldsymbol{\cdot}C$.
Thus $A\boldsymbol{\cdot}B\boldsymbol{\cdot}C$ might be any one among the six matrices: $A\circ B\circ C$, $ABC$, $A\circ(BC)$, $A(B\circ C)$, $(A\circ B)C$ or $(AB)\circ C$.

For any $A\in\M_n(\IC)$, the singular values of $A$ are the nonnegative square roots of the eigenvalues of $A^*A$ (or $AA^*$ as for any two $A, B$, $AB$ and $BA$ have the same eigenvalues). In other words, singular values of $A$ are precisely the eigenvalues of $|A|$ counted with multiplicities. Here $|A|$ denotes the unique positive square root of the matrix $A^*A$. The singular values of $A\in\M_n(\IC)$ are denoted by $\sigma_i(A)$, $1\leq i\leq n$ and they are always assumed to be arranged in decreasing order, i.e.
\be\label{P6eq10}
\sigma_1(A)\geq\sigma_2(A)\geq\cdots\geq\sigma_n(A)\geq 0.
\ee
For any $A\in\M_n(\IC)$ (see \cite[p. 291]{Horn1} and \cite[p. 421, Ex. 4]{Horn1}),
\be\label{P6eq11}
\sum_{i, j=1}^n|A_{ij}|^2
=\sum_{i=1}^n\sigma_i^2(A)
=\sum_{i=1}^n\sigma_i(A^*A)
=\mbox{tr}(A^*A).
\ee
For any $A\in\M_n(\IC)$, there exists a unitary $U\in\M_n(\IC)$ such that
\be\label{P6eq12}
A=U|A|
\ee
(see \cite[p. 413, Corollary 7.3.3]{Horn1}). Also, there exist unitary $V, W\in\M_n(\IC)$ such that
\be\label{P6eq13}
A=V\Sigma W^*
\ee
(cf. \cite[p. 414, Theorem 7.3.5]{Horn1}), where $\Sigma$ is a diagonal matrix with the diagonal entries $\Sigma_{ii}=\sigma_i(A)$, $1\leq i\leq n$.
A norm $\|.\|$ on $\M_n(\IC)$ is said to be \emph{unitarily invariant} if $\|UAV\|=\|A\|$ for any unitary $U, V$ in $\M_n(\IC)$, and for all $A\in\M_n(\IC)$. Clearly, according to the above notations, $\|A\|=\|\Sigma\|$. At this point, we record the following useful results.

\begin{theorem}\label{TheoB}\cite[Theorems 1.4, 1.2]{Li}
For any $A, B, C\in\M_n(\IC)$ and any unitarily invariant norm $\|.\|$ on $\M_n(\IC)$,
$$
\|A\|^2\leq\|B\|\|C\|\iff
\left(\sum_{i=1}^k\sigma_i(A)\right)^2\leq
\left(\sum_{i=1}^k\sigma_i(B)\right)\left(\sum_{i=1}^k\sigma_i(C)\right)\,\,\forall\, 1\leq k\leq n.
$$
In particular, taking $B=C$ we re-obtain the well known result:
$$
\|A\|\leq\|B\|\iff
\left(\sum_{i=1}^k\sigma_i(A)\right)\leq
\left(\sum_{i=1}^k\sigma_i(B)\right)\,\,\forall\, 1\leq k\leq n.
$$
\end{theorem}
\begin{theorem}\label{TheoC}\cite[p. 206, Corollary 3.5.10]{Horn}
Let $A, B\in\M_n(\IC)$ and $\|.\|$ be any unitarily invariant norm in $\M_n(\IC)$. Also, let $E_{11}\in\M_n(\IC)$ have the entry $1$ in position $(1,1)$ and zeroes elsewhere. Then
$\left(a\right)\|AB^*\|\leq\sigma_1(A)\|B\|$ and
$\left(b\right)\|A\|\geq\sigma_1(A)\|E_{11}\|$.
\end{theorem}

A function $\phi:\IC^n\to[0,\infty)$ is said to be a \emph{symmetric gauge function} if it is a norm on $\IC^n$, satisfying the following two additional properties:
\bee
\item
$\phi$ is an absolute norm, i.e. $\phi(x_1, x_2, x_3,\cdots, x_n)=\phi(|x_1|, |x_2|, |x_3|,\cdots, |x_n|)$ for any $(x_1, x_2, x_3,\cdots, x_n)\in\IC^n$.
\item
$\phi(Px)=\phi(x)$ for all $x\in\IC^n$ and for any permutation matrix $P\in\M_n(\IC)$. In other words, $\phi(x)=\phi(y)$ for any $x=(x_1, x_2, x_3,\cdots, x_n)\in\IC^n$ and any rearrangement $y=(y_1, y_2, y_3,\cdots, y_n)$ of
$(x_1, x_2, x_3,\cdots, x_n)$.
\eee
The following result would be of use later in this article.

\begin{theorem}\label{TheoD} \cite[Theorem 7.4.45 \& (7.4.46)]{Horn1}
Suppose $x=(x_1, x_2, x_3,\cdots, x_n)\in\IC^n$ and
$y=(y_1, y_2, y_3,\cdots, y_n)\in\IC^n$ are such that
$|x_1|\leq|x_2|\leq|x_3|\leq\cdots\leq|x_n|$ and
$|y_1|\leq|y_2|\leq|y_3|\leq\cdots\leq|y_n|$. Let $\phi$ be any symmetric gauge function on $\IC^n$. Then
$$
\phi(x)\leq\phi(y)\iff
\sum_{i=t}^{n}|x_i|\leq\sum_{i=t}^{n}|y_i|\,\, \forall\,\,
1\leq t\leq n.
$$
\end{theorem}
Since $\M_n(\IC)$ can be identified with $\IC^{n^2}$ as a complex vector space, $\phi(A)$ makes sense for any $A\in\M_n(\IC)$ if $A$ is represented as an element of $\IC^{n^2}$. In fact, we can write $A=(\alpha_1, \alpha_2, \cdots, \alpha_{n^2})$, where \lq$\alpha_i$\rq s are precisely the entries of matrix $A$ written in any order; and then $\phi(A):=\phi(\alpha_1, \alpha_2, \cdots, \alpha_{n^2})$.
The reader is urged to glance through \cite{Horn1, Horn} for a detailed study of unitarily invariant norms and symmetric gauge functions.

We end this section by introducing another interesting real valued function defined on $\M_n(\IC)$ that is not a norm.
Let $S_n$ denote the symmetric group of degree $n$ and $H$ be a subgroup of $S_n$ of order $h$. Let $\chi$ be a character of degree 1 on $H$,
i.e. a nontrivial homomorphism of $H$ into the complex numbers.
For any $A\in\M_n(\IC)$, we define the \emph{generalized matrix function} $M_{\chi}$ (cf. \cite{Marc}) by
$$
M_{\chi}(A)= \sum_{\sigma\in H}\chi(\sigma)\prod_{i=1}^n A_{i\sigma(i)}.
$$
We will make use of the following result in our proofs.

\begin{theorem}\label{TheoE}\cite[Corollary 3.2]{Marc}
$|M_\chi(A)|^2\leq(1/n)\sum_{i=1}^n\sigma_i^{2n}(A)$
for any $A\in\M_n(\IC)$.
\end{theorem}
\section{Main results and their proofs}
We need the following lemma for obtaining the subsequent results. This lemma might appear somewhat straightforward from already known facts, but we include a proof for the sake of completeness.
\blem\label{P6lem1}
For any $A, B, C\in\M_n(\IC)$, $n\in\IN$ and for any $k=1, 2, \cdots n$ we have
\be\label{P6eq1}
\sum_{i=1}^k\sigma_i^m(A\boldsymbol{\cdot}B\boldsymbol{\cdot}C)
\leq\sum_{i=1}^k\sigma_i^m(A)\sigma_i^m(B)\sigma_i^m(C)
\ee
for any $m\in\IN$.
\elem
\bpf
We give a detailed proof for the case
$A\boldsymbol{\cdot}B\boldsymbol{\cdot}C=
A\boldsymbol{\cdot}(B\boldsymbol{\cdot}C)$ only, as the proof for the case
$A\boldsymbol{\cdot}B\boldsymbol{\cdot}C=
(A\boldsymbol{\cdot}B)\boldsymbol{\cdot}C$
is exactly similar.
From \cite[Theorem 3.3.14(a), Theorem 5.5.4]{Horn} it is known that for any $B, C\in\M_n(\IC)$ and for all $k$ satisfying $1\leq k\leq n$,
$$
\sum_{i=1}^k\sigma_i(B\boldsymbol{\cdot}C)
\leq\sum_{i=1}^k\sigma_i(B)\sigma_i(C).
$$
Since $\sigma_i(A)$ is a monotonically decreasing finite sequence of nonnegative real numbers over $i$, using summation by parts combined with the above inequality we have
\be\label{P6eq20}
\sum_{i=1}^k\sigma_i(A\boldsymbol{\cdot}B\boldsymbol{\cdot}C)
\leq\sum_{i=1}^k\sigma_i(A)\sigma_i(B\boldsymbol{\cdot}C)
\leq\sum_{i=1}^k\sigma_i(A)\sigma_i(B)\sigma_i(C).
\ee
In other words, $(\ref{P6eq1})$ remains true for $m=1$ and for all $k$ satisfying $1\leq k\leq n$.
Let $(\ref{P6eq1})$ be true for $m=m_1\in\IN$ and for all $k$, $1\leq k\leq n$, i.e.
\be\label{P6eq21}
\sum_{i=1}^k\sigma_i^{m_1}(A\boldsymbol{\cdot}B\boldsymbol{\cdot}C)
\leq\sum_{i=1}^k\sigma_i^{m_1}(A)\sigma_i^{m_1}(B)\sigma_i^{m_1}(C).
\ee
Now both $\sigma_i(A\boldsymbol{\cdot}B\boldsymbol{\cdot}C)$ and $\sigma_i^{m_1}(A)\sigma_i^{m_1}(B)\sigma_i^{m_1}(C)$ are monotonically decreasing finite sequences of nonnegative real numbers over $i$. Therefore, using the summation by parts with $(\ref{P6eq21})$, and then again with $(\ref{P6eq20})$ we get
\begin{align*}
&\sum_{i=1}^k\sigma_i^{m_1+1}(A\boldsymbol{\cdot}B\boldsymbol{\cdot}C)\\
&=\sum_{i=1}^k\sigma_i^{m_1}(A\boldsymbol{\cdot}B\boldsymbol{\cdot}C)
\sigma_i(A\boldsymbol{\cdot}B\boldsymbol{\cdot}C)\\
&\leq \sum_{i=1}^k\sigma_i^{m_1}(A)\sigma_i^{m_1}(B)\sigma_i^{m_1}(C)\sigma_i(A\boldsymbol{\cdot}B\boldsymbol{\cdot}C)\\
&\leq\sum_{i=1}^k\sigma_i^{m_1}(A)\sigma_i^{m_1}(B)\sigma_i^{m_1}(C)\sigma_i(A)\sigma_i(B)\sigma_i(C)\\
&=\sum_{i=1}^k\sigma_i^{m_1+1}(A)\sigma_i^{m_1+1}(B)\sigma_i^{m_1+1}(C)
\end{align*}
for any $k$ satisfying $1\leq k\leq n$. The proof is now complete by the principle of mathematical induction.
\epf

We now state and prove the first theorem of this article. It should be mentioned here that for any complex matrix $A$, $\ov{A}$ denotes the conjugate of $A$, i.e. the $(i,j)$-th entry of $\ov{A}$ is precisely the complex conjugate of the $(i,j)$-th entry of $A$.
\bthm\label{P6thm1}
Let $f\in L^1(G, \IC)$ and $\emph{Re}(f)\leq 1$ a.e. on $G$, $0\leq\widehat{f}(\pi_0)<1$.
Also let $\{R_n\}_{n\in\Lambda}$ be a sequence of matrices such that $R_n\in\M_{d_{\pi_n}}(\IC)$. Then for all $x\in G$, we have the following:
\bee
\item[(i)]~If $\|.\|_{\pi_n}$ is any unitarily invariant norm on $\M_{d_{\pi_n}}(\IC)$, then for $\sum_{n\in\Lambda} \|R_n\|_{\pi_n}\leq 1/2$ we have
\be\label{P6eq29}
\widehat{f}(\pi_0)+
\sum_{n\in\Lambda}\left\|\left(\widehat{f}(\pi_n)+\ov{\widehat{f}(\ov{\pi_n})}\right)
\boldsymbol{\cdot}\pi_n(x)\boldsymbol{\cdot}R_{n}\right\|_{\pi_n}\leq 1.
\ee
\item[(ii)]~If $\phi_n$ is any symmetric gauge function on $\IC^{d^2_{\pi_n}}$, then for $\sum_{n\in\Lambda} d_{\pi_n}^2\phi_n(R_n)\leq 1/2$ we have
\be\label{P6eq30}
\widehat{f}(\pi_0)+\sum_{n\in\Lambda}\phi_n\left(\left(\widehat{f}(\pi_n)+\ov{\widehat{f}(\ov{\pi_n})}\right)
\boldsymbol{\cdot}\pi_n(x)\boldsymbol{\cdot}R_{n}\right)\leq 1.
\ee
\item[(iii)]~
If $M_{\chi_n}$ is a generalized matrix function on $\M_{d_{\pi_n}}(\IC)$, where
$\chi_n$ is a character of degree $1$ on a subgroup $H$ of $S_{d_{\pi_n}}$, then for
$$
\sum_{n\in\Lambda}\frac{\|R_n\|_{\pi_n}}{\|E^{(n)}_{11}\|_{\pi_n}}\leq 1/2
$$
we have
\be\label{P6eq31}
\widehat{f}(\pi_0)+\sum_{n\in\Lambda}
\left|M_{\chi_n}\left(\left(\widehat{f}(\pi_n)+\ov{\widehat{f}(\ov{\pi_n})}\right)
\boldsymbol{\cdot}\pi_n(x)\boldsymbol{\cdot}R_{n}\right)\right|^{\frac{1}{d_{\pi_n}}}\leq 1.
\ee
\eee
If $\widehat{f}(\pi_0)=1$ then equality holds at $(\ref{P6eq29})$, $(\ref{P6eq30})$ and $(\ref{P6eq31})$ for any sequence $\{R_n\}_{n\in\Lambda}$, i.e. \lq$R_n$\rq s need not satisfy any inequality as in the parts $(i)$, $(ii)$ and $(iii)$ above.
\ethm
The point of departure in the following proof of Theorem \ref{P6thm1} is an extension of the technique used for establishing \cite[Theorem 2.1]{Paul1} (or \cite[Theorem 4.1]{Paul3}).
\bpf[Proof of Theorem \ref{P6thm1}]
We first complete the proofs of all three parts for $0\leq\widehat{f}(\pi_0)<1$.
Using the expression at $(\ref{P6eq8})$, the $(i,j)$-th element of $\ov{\widehat{f}(\ov{\pi_n})}$ is
$$
\ov{\widehat{f}(\ov{\pi_n})}_{ij}
=\ov{\int_G f(x)\pi_n(x)_{ji}\,dx}
=\int_G \ov{f(x)}\,\ov{\pi_n(x)_{ji}}\,dx
$$
for any $n\in\Lambda\cup\{0\}$. Hence, according to $(\ref{P6eq6})$, we have
\be\label{P6eq22}
\ov{\widehat{f}(\ov{\pi_n})}
=\widehat{\ov{f}}(\pi_n)
=\int_G\ov{f(x)}\pi_n^*(x)\,dx.
\ee
Since $\pi_0(x)=1$ for all $x\in G$, by virtue of the Peter-Weyl Theorem (cf. \cite[p. 133]{Foll}) it is immediate that
$$
\int_G\pi_n^*(x)\,dx=0
$$
for any $n\in\Lambda$.
Therefore, adding $(\ref{P6eq6})$ and $(\ref{P6eq22})$ we have
$$
-\left(\widehat{f}(\pi_n)+\ov{\widehat{f}(\ov{\pi_n})}\right)
=2\int_G\mbox{Re}\left(1-f(x)\right)\pi_n^*(x)\,dx.
$$
As $-\left(\widehat{f}(\pi_n)+\ov{\widehat{f}(\ov{\pi_n})}\right)\in\M_{d_{\pi_n}}(\IC)$,
there exists a unitary $U_n\in\M_{d_{\pi_n}}(\IC)$ (see ($\ref{P6eq12}$)) such that
$$
-\left(\widehat{f}(\pi_n)+\ov{\widehat{f}(\ov{\pi_n})}\right)
=U_n\left|\widehat{f}(\pi_n)+\ov{\widehat{f}(\ov{\pi_n})}\right|.
$$
Hence
\be\label{P6eq23}
\left\langle\left|\widehat{f}(\pi_n)+\ov{\widehat{f}(\ov{\pi_n})}\right|u, u\right\rangle
=2\int_G\mbox{Re}\left(1-f(x)\right)\left\langle U_n^*\pi_n^*(x)u, u\right\rangle\,dx
\ee
for any $u\in\IH_{\pi_n}$, $\langle . , .\rangle$ being the inner product on $\IH_{\pi_n}$. Since $U_n^*\pi_n^*(x)$ is unitary for all $x\in G$,
$\|U_n^*\pi_n^*(x)u\|_h=\|u\|_h$
for any $u\in\IH_{\pi_n}$, which further implies
$$|\langle U_n^*\pi_n^*(x)u, u\rangle|\leq \|U_n^*\pi_n^*(x)u\|_h\|u\|_h=\|u\|_h^2.
$$
Here $\|.\|_h$ denotes the norm on the Hilbert space $\IH_{\pi_n}$ induced by the inner product $\langle . , . \rangle$.
Thus from $(\ref{P6eq23})$ we get
\begin{align*}
&\left\langle\left|\widehat{f}(\pi_n)+\ov{\widehat{f}(\ov{\pi_n})}\right|u, u\right\rangle\\
&\leq2\int_G\mbox{Re}\left(1-f(x)\right)\left|\left\langle U_n^*\pi_n^*(x)u, u\right\rangle\right|\,dx\\
&\leq \left(2\int_G\mbox{Re}(1-f(x))\,dx\right)\left\langle u, u\right\rangle\\
&=2\left(1-\widehat{f}(\pi_0)\right)\left\langle u, u\right\rangle.
\end{align*}
Now for any eigenvalue $\sigma_i\left(\widehat{f}(\pi_n)+\ov{\widehat{f}(\ov{\pi_n})}\right)$ of $\left|\widehat{f}(\pi_n)+\ov{\widehat{f}(\ov{\pi_n})}\right|$,
let $u$ be the corresponding nonzero eigenvector. It is easily seen from the above inequality that
$\sigma_i\left(\widehat{f}(\pi_n)+\ov{\widehat{f}(\ov{\pi_n})}\right)\left\langle u, u\right\rangle\leq 2\left(1-\widehat{f}(\pi_0)\right)
\left\langle u, u\right\rangle,
$
i.e.
\be\label{P6eq2}
\sigma_i\left(\widehat{f}(\pi_n)+\ov{\widehat{f}(\ov{\pi_n})}\right)
\leq 2\left(1-\widehat{f}(\pi_0)\right)
\ee
for any $i\in\IN$ such that $1\leq i\leq d_{\pi_n}$.
Now using the Lemma \ref{P6lem1}, the fact that $\sigma_i(\pi_n(x))=1$ for any $i\in\IN$ satisfying $1\leq i\leq d_{\pi_n}$, and $(\ref{P6eq2})$ we have
\begin{align*}
&\sum_{i=1}^k\sigma_i\left(\left(\widehat{f}(\pi_n)+\ov{\widehat{f}(\ov{\pi_n})}\right)
\boldsymbol{\cdot}\pi_n(x)\boldsymbol{\cdot}R_{n}\right)\\
&\leq\sum_{i=1}^k\sigma_i\left(\widehat{f}(\pi_n)+\ov{\widehat{f}(\ov{\pi_n})}\right)
\sigma_i\left(\pi_n(x)\right)\sigma_i\left(R_{n}\right)\\
&\leq 2\left(1-\widehat{f}(\pi_0)\right)\sum_{i=1}^k\sigma_i\left(R_n\right)
\end{align*}
for any $k\in\IN$ satisfying $1\leq k\leq d_{\pi_n}$. From Theorem \ref{TheoB}, we conclude that
$$
\left\|\left(\widehat{f}(\pi_n)+\ov{\widehat{f}(\ov{\pi_n})}\right)
\boldsymbol{\cdot}\pi_n(x)\boldsymbol{\cdot}R_{n}\right\|_{\pi_n}
\leq 2\left(1-\widehat{f}(\pi_0)\right)\left\|R_n\right\|_{\pi_n}.
$$
Using the above inequality, part (i) of Theorem \ref{P6thm1} follows from a simple computation.

We now concentrate on part (ii). Suppose that the representations of
$\left(\widehat{f}(\pi_n)+\ov{\widehat{f}(\ov{\pi_n})}\right)
\boldsymbol{\cdot}\pi_n(x)\boldsymbol{\cdot}R_{n}$ and $R_n$
as elements of $\IC^{d_{\pi_n}^2}$ are as follows:
$$
\left(\widehat{f}(\pi_n)+\ov{\widehat{f}(\ov{\pi_n})}\right)
\boldsymbol{\cdot}\pi_n(x)\boldsymbol{\cdot}R_{n}=\left(a_1, a_2, \cdots, a_{d_{\pi_n}^2}\right)\, \mbox{and}\,
R_n=\left(b_1, b_2, \cdots, b_{d_{\pi_n}^2}\right).
$$
Now applying the permutation invariance of $\phi_n$, the Cauchy-Schwarz inequality, $(\ref{P6eq11})$, Lemma \ref{P6lem1}, $(\ref{P6eq2})$ and the fact that $\sigma_i(\pi_n(x))=1$, and then again $(\ref{P6eq11})$ appropriately, we have
\begin{align*}
&\phi_n\left(\left(\widehat{f}(\pi_n)+\ov{\widehat{f}(\ov{\pi_n})}\right)
\boldsymbol{\cdot}\pi_n(x)\boldsymbol{\cdot}R_{n}\right)=\phi_n\left(a_1, a_2, \cdots, a_{d_{\pi_n}^2}\right)\\
&=\phi_n\left(a_1(1, 0, 0, \cdots, 0)+ a_2(0, 1, 0, \cdots, 0)+\cdots+a_{d_{\pi_n}^2}(0, 0, 0, \cdots, 1)\right)\\
&\leq\phi_n\left(1, 0, 0, \cdots, 0\right)\sum_{i=1}^{d_{\pi_n}^2}\left|a_i\right|\\
&\leq d_{\pi_n}\phi_n(1, 0, 0, \cdots, 0)\sqrt{\sum_{i=1}^{d_{\pi_n}^2}|a_i|^2}\\
&=d_{\pi_n}\phi_n\left(1, 0, 0, \cdots, 0\right)\sqrt{\sum_{i=1}^{d_{\pi_n}}\sigma_i^2\left(\left(\widehat{f}(\pi_n)+\ov{\widehat{f}(\ov{\pi_n})}\right)\boldsymbol{\cdot}\pi_n(x)\boldsymbol{\cdot}R_n\right)}\\
&\leq d_{\pi_n}\phi_n\left(1, 0, 0, \cdots, 0\right)\sqrt{\sum_{i=1}^{d_{\pi_n}}\sigma_i^2\left(\widehat{f}(\pi_n)+\ov{\widehat{f}(\ov{\pi_n})}\right)\sigma_i^2\left(\pi_n(x)\right)\sigma_i^2\left(R_n\right)}\\
&\leq 2\left(1-\widehat{f}(\pi_0)\right)d_{\pi_n}\phi_n\left(1, 0, 0, \cdots, 0\right)\sqrt{\sum_{i=1}^{d_{\pi_n}}\sigma_i^2\left(R_n\right)}\\
&=2\left(1-\widehat{f}(\pi_0)\right)d_{\pi_n}\phi_n\left(1, 0, 0, \cdots, 0\right)\sqrt{\sum_{i=1}^{d_{\pi_n}^2}\left|b_i\right|^2}\\
&\leq 2\left(1-\widehat{f}(\pi_0)\right)d_{\pi_n}^2\phi_n\left(1, 0, 0, \cdots, 0\right)\left|b_{\mbox{max}}\right|,
\end{align*}
where $\left|b_{\mbox{max}}\right|=\mbox{max}_{1\leq i\leq d_{\pi_n}^2}\left|b_i\right|$.
Now we rearrange the entries of $\left(b_1, b_2, \cdots, b_{d_{\pi_n}^2}\right)$ to get the new vector
$\left(b_{r_1}, b_{r_2},\cdots, b_{r_{d_{\pi_n}^2}}\right)$. The rearrangement is done in such a way that the finite sequence $\left\{\left|b_{r_i}\right|\right\}_{i=1}^{d_{\pi_n}^2}$ is monotonically increasing. It is immediately seen that
for any $t=1, 2, \cdots, d_{\pi_n}^2$,
$$
\left|b_{\mbox{max}}\right|=\left|b_{r_{d_{\pi_n}^2}}\right|\leq\sum_{i=t}^{d_{\pi_n}^2}\left|b_{r_i}\right|.
$$
Using Theorem \ref{TheoD}, we observe that
$$
\phi_n\left(0, 0, \cdots, \left|b_{\mbox{max}}\right|\right)\leq\phi_n\left(b_{r_1}, b_{r_2}, \cdots, b_{r_{d_{\pi_n}^2}}\right),
$$
or equivalently, using the permutation invariance of $\phi_n$:
\be\label{P6eq5}
\phi_n\left(1, 0, 0, \cdots, 0\right)\left|b_{\mbox{max}}\right|\leq\phi_n\left(b_1, b_2, \cdots, b_{d_{\pi_n}^2}\right)=
\phi_n\left(R_n\right).
\ee
Therefore
$$
\phi_n\left(\left(\widehat{f}(\pi_n)+\ov{\widehat{f}(\ov{\pi_n})}\right)
\boldsymbol{\cdot}\pi_n(x)\boldsymbol{\cdot}R_{n}\right)
\leq2\left(1-\widehat{f}(\pi_0)\right)d_{\pi_n}^2\phi_n\left(R_n\right).
$$
The remaining portion of the proof can be completed by a direct calculation using the above inequality.

For part (iii), we need to use Theorem \ref{TheoE}, Lemma \ref{P6lem1}, $\sigma_i(\pi_n(x))=1$ and $(\ref{P6eq2})$, $(\ref{P6eq10})$ and finally Theorem \ref{TheoC}(b) one by one. As a result, we get
\begin{align*}
&\left|M_{\chi_n}\left(\left(\widehat{f}(\pi_n)+\ov{\widehat{f}(\ov{\pi_n})}\right)
\boldsymbol{\cdot}\pi_n(x)\boldsymbol{\cdot}R_{n}\right)\right|^2\\
&\leq \frac{1}{d_{\pi_n}}\sum_{i=1}^{d_{\pi_n}}
\sigma_i^{2d_{\pi_n}}\left(\left(\widehat{f}(\pi_n)+\ov{\widehat{f}(\ov{\pi_n})}\right)
\boldsymbol{\cdot}\pi_n(x)\boldsymbol{\cdot}R_{n}\right)
\\
&\leq \frac{1}{d_{\pi_n}}\sum_{i=1}^{d_{\pi_n}}
\sigma_i^{2d_{\pi_n}}\left(\widehat{f}(\pi_n)+\ov{\widehat{f}(\ov{\pi_n})}\right)
\sigma_i^{2d_{\pi_n}}\left(\pi_n(x)\right)\sigma_i^{2d_{\pi_n}}\left(R_n\right)\\
&\leq \frac{\left(2\left(1-\widehat{f}(\pi_0)\right)\right)^{2d_{\pi_n}}}{d_{\pi_n}}
\sum_{i=1}^{d_{\pi_n}}\sigma_i^{2d_{\pi_n}}\left(R_n\right)\\
&\leq \left(2\left(1-\widehat{f}(\pi_0)\right)\right)^{2d_{\pi_n}}\sigma_1^{2d_{\pi_n}}\left(R_n\right)\\
&\leq \left(2\left(1-\widehat{f}(\pi_0)\right)\right)^{2d_{\pi_n}}
\frac{\|R_n\|_{\pi_n}^{2d_{\pi_n}}}{\|E_{11}^{(n)}\|_{\pi_n}^{2d_{\pi_n}}}\,.
\end{align*}
In other words,
$$
\left|M_{\chi_n}\left(\left(\widehat{f}(\pi_n)+\ov{\widehat{f}(\ov{\pi_n})}\right)
\boldsymbol{\cdot}\pi_n(x)\boldsymbol{\cdot}R_{n}\right)\right|^{\frac{1}{d_{\pi_n}}}
\leq 2\left(1-\widehat{f}(\pi_0)\right)
\frac{\|R_n\|_{\pi_n}}{\|E_{11}^{(n)}\|_{\pi_n}},
$$
from which the proof now follows by a simple calculation.

Finally, if $\widehat{f}(\pi_0)=1$ then $(\ref{P6eq2})$ implies that $\sigma_i\left(\widehat{f}(\pi_n)+\ov{\widehat{f}(\ov{\pi_n})}\right)=0$ for all $1\leq i\leq d_{\pi_n}$, and hence $\widehat{f}(\pi_n)+\ov{\widehat{f}(\ov{\pi_n})}=0$ (here $0$ means the zero matrix in the vector space $\M_{d_{\pi_n}}(\IC)$) for any $n\in\Lambda$. Thus equality occurs at $(\ref{P6eq29})$, $(\ref{P6eq30})$ and $(\ref{P6eq31})$, irrespective of the choices of $R_n$. This finishes the proof of Theorem \ref{P6thm1}.
\epf

We now proceed to prove the next theorem of this article.
\bthm\label{P6thm2}
Suppose $f\in L^2(G, \IC)$, $\|f\|_{L^2(G, \IC)}\leq 1$ and $0\leq\widehat{f}(\pi_0)<1$. Further suppose that $R_n, \|.\|_{\pi_n}, \phi_n$ and $M_{\chi_n}$ bear the same meaning as in the statement of Theorem~\ref{P6thm1}. Then for all $x\in G$, the following are true:
\bee
\item[(i)]~ If $(1+\widehat{f}(\pi_0))\left(1+\sum_{n\in\Lambda} d_{\pi_n}\|R_n\|_{\pi_n}^2\right)\leq 2$, then
\be\label{P6eq32}
\widehat{f}(\pi_0)+\sum_{n\in\Lambda} d_{\pi_n}\|\widehat{f}(\pi_n)\boldsymbol{\cdot}\pi_n(x)\boldsymbol{\cdot}R_n\|_{\pi_n}\leq 1.
\ee
\item[(ii)]~If $(1+\widehat{f}(\pi_0))\left(1+\sum_{n\in\Lambda} d^5_{\pi_n}\phi_n^2(R_n)\right)\leq 2$, then
\be\label{P6eq33}
\widehat{f}(\pi_0)+\sum_{n\in\Lambda} d_{\pi_n}\phi_n(\widehat{f}(\pi_n)\boldsymbol{\cdot}\pi_n(x)\boldsymbol{\cdot}R_n)\leq 1.
\ee
\item[(iii)]~ If
$$
\left(1+\widehat{f}(\pi_0)\right)\left(1+\sum_{n\in\Lambda} d_{\pi_n}\frac{\|R_n\|_{\pi_n}^2}{\|E^{(n)}_{11}\|_{\pi_n}^2}\right)\leq 2,
$$
then
\be\label{P6eq34}
\widehat{f}(\pi_0)+\sum_{n\in\Lambda}
d_{\pi_n}|M_{\chi_n}(\widehat{f}(\pi_n)\boldsymbol{\cdot}\pi_n(x)\boldsymbol{\cdot}R_n)|^{\frac{1}{d_{\pi_n}}}\leq 1.
\ee
\eee
If $\widehat{f}(\pi_0)=1$ then equality holds at $(\ref{P6eq32})$, $(\ref{P6eq33})$ and $(\ref{P6eq34})$ for any sequence $\{R_n\}_{n\in\Lambda}$, i.e. \lq$R_n$\rq s need not satisfy any inequality as in the parts $(i)$, $(ii)$ and $(iii)$ above.
\ethm
\bpf[Proof of Theorem \ref{P6thm2}]
We first complete the proofs of all three parts under the assumption $0\leq\widehat{f}(\pi_0)<1$.
To prove part (i), we start with the fact that for any unitarily invariant norm $\|.\|_{\pi_n}$,
$$
\|\widehat{f}(\pi_n)^*\widehat{f}(\pi_n)\|_{\pi_n}=
\|\Sigma_n\|_{\pi_n}
\leq\sum_{i=1}^{d_{\pi_n}}
\sigma_i(\widehat{f}(\pi_n)^*\widehat{f}(\pi_n))\|E_{ii}^{(n)}\|_{\pi_n},
$$
where $\Sigma_n\in\M_{d_{\pi_n}}(\IC)$ is a diagonal matrix with $\sigma_i(\widehat{f}(\pi_n)^*\widehat{f}(\pi_n))$ as the $(i,i)$-th diagonal entry, $1\leq i\leq d_{\pi_n}$. Further, it is evident from Theorem \ref{TheoB} that
$\|E_{11}^{(n)}\|_{\pi_n}\leq\|E_{ii}^{(n)}\|_{\pi_n}$ and
$\|E_{ii}^{(n)}\|_{\pi_n}\leq\|E_{11}^{(n)}\|_{\pi_n}$, as all $E_{ii}^{(n)}$ has only one singular value $1$ and other singular values $0$. Therefore,
$\|E_{11}^{(n)}\|_{\pi_n}=\|E_{ii}^{(n)}\|_{\pi_n}$ for any $i$. Using this and $(\ref{P6eq11})$, we get
\be\label{P6eq3}
\|\widehat{f}(\pi_n)^*\widehat{f}(\pi_n)\|_{\pi_n}
\leq\mbox{tr}(\widehat{f}(\pi_n)^*\widehat{f}(\pi_n))\|E_{11}^{(n)}\|_{\pi_n}.
\ee
Making use of Lemma \ref{P6lem1}, the fact that $\sigma_i(\pi_n(x))=1$ and the Cauchy-Schwarz inequality, we observe that
\begin{align*}
&\left(\sum_{i=1}^k\sigma_i(\widehat{f}(\pi_n)\boldsymbol{\cdot}\pi_n(x)\boldsymbol{\cdot}R_n)\right)^2\\
&\leq\left(\sum_{i=1}^k\sigma_i(\widehat{f}(\pi_n))\sigma_i(\pi_n(x))\sigma_i(R_n)\right)^2\\
&=\left(\sum_{i=1}^k\sigma_i(\widehat{f}(\pi_n))\sigma_i(R_n)\right)^2\\
&\leq\left(\sum_{i=1}^k\sigma_i^2(\widehat{f}(\pi_n))\right)\left(\sum_{i=1}^k\sigma_i^2(R_n)\right)\\
&=\left(\sum_{i=1}^k\sigma_i(\widehat{f}(\pi_n)^*\widehat{f}(\pi_n))\right)\left(\sum_{i=1}^k\sigma_i(R_nR_n^*)\right)
\end{align*}
for any $k\in\IN$ satisfying $1\leq k\leq d_{\pi_n}$.
It follows from Theorem \ref{TheoB} that
\be\label{P6eq4}
\|\widehat{f}(\pi_n)\boldsymbol{\cdot}\pi_n(x)\boldsymbol{\cdot}R_n\|_{\pi_n}^2
\leq\|\widehat{f}(\pi_n)^*\widehat{f}(\pi_n)\|_{\pi_n}\|R_nR_n^*\|_{\pi_n}.
\ee
A use of Theorem \ref{TheoC} with the fact that $\|R_n^*\|_{\pi_n}=\|R_n\|_{\pi_n}$ gives
$$
\|R_nR_n^*\|_{\pi_n}\leq\sigma_1(R_n)\|R_n\|_{\pi_n}\leq\frac{\|R_n\|_{\pi_n}^2}{\|E_{11}^{(n)}\|_{\pi_n}}.
$$
Combining the above inequality with $(\ref{P6eq3})$ and $(\ref{P6eq4})$, we obtain
$$
\|\widehat{f}(\pi_n)\boldsymbol{\cdot}\pi_n(x)\boldsymbol{\cdot}R_n\|_{\pi_n}
\leq\sqrt{\mbox{tr}(\widehat{f}(\pi_n)^*\widehat{f}(\pi_n))}\|R_n\|_{\pi_n}.
$$
From $(\ref{P6eq9})$ we have
\be\label{P6eq24}
\|f\|_{L^2(G, \IC)}^2\leq 1
\iff\widehat{f}(\pi_0)^2+\sum_{n\in\Lambda}  d_{\pi_n}\mbox{tr}(\widehat{f}(\pi_n)^*\widehat{f}(\pi_n))
\leq 1.
\ee
Now, using the above two inequalities along with the Cauchy-Schwarz inequality, we get
\begin{align*}
&\sum_{n\in\Lambda} d_{\pi_n}\|\widehat{f}(\pi_n)\boldsymbol{\cdot}\pi_n(x)\boldsymbol{\cdot}R_n\|_{\pi_n}\\
&\leq\sum_{n\in\Lambda}\sqrt{d_{\pi_n}\mbox{tr}(\widehat{f}(\pi_n)^*\widehat{f}(\pi_n))}\sqrt{d_{\pi_n}}\|R_n\|_{\pi_n}\\
&\leq\sqrt{\sum_{n\in\Lambda} d_{\pi_n}\mbox{tr}(\widehat{f}(\pi_n)^*\widehat{f}(\pi_n))}\sqrt{\sum_{n\in\Lambda} d_{\pi_n}\|R_n\|_{\pi_n}^2}\\
&\leq\sqrt{1-\widehat{f}(\pi_0)^2}\sqrt{\sum_{n\in\Lambda} d_{\pi_n}\|R_n\|_{\pi_n}^2}.
\end{align*}
Rest of the proof is a straightforward calculation using the above inequality.

Most of the work required for proving part (ii) of this theorem is already done in the proof of Theorem \ref{P6thm1}. We only have to replace $\widehat{f}(\pi_n)+\ov{\widehat{f}(\ov{\pi_n})}$ by $\widehat{f}(\pi_n)$ and then adopt the similar lines of arguments as in the proof of part (ii) of Theorem \ref{P6thm1} (with necessary little changes). With an application of $(\ref{P6eq11})$, it follows that
\begin{align*}
&\phi_n(\widehat{f}(\pi_n)\boldsymbol{\cdot}\pi_n(x)\boldsymbol{\cdot}R_n)\\
&\leq d_{\pi_n}\phi_n(1, 0, 0, \cdots, 0)\sqrt{\sum_{i=1}^{d_{\pi_n}}\sigma_i^2(\widehat{f}(\pi_n))\sigma_i^2(\pi_n(x))\sigma_i^2(R_n)}\\
&= d_{\pi_n}\phi_n(1, 0, 0, \cdots, 0)\sqrt{\sum_{i=1}^{d_{\pi_n}}\sigma_i^2(\widehat{f}(\pi_n))\sigma_i^2(R_n)}\\
&\leq d_{\pi_n}\phi_n(1, 0, 0, \cdots, 0)\sqrt{\sum_{i=1}^{d_{\pi_n}}\sigma_i^2(\widehat{f}(\pi_n))}\sqrt{\sum_{i=1}^{d_{\pi_n}}\sigma_i^2(R_n)}\\
&\leq d_{\pi_n}^2\sqrt{\mbox{tr}(\widehat{f}(\pi_n)^*\widehat{f}(\pi_n))}\phi_n(1, 0, 0, \cdots, 0)|b_{\mbox{max}}|,
\end{align*}
$|b_{\mbox{max}}|$ as defined in the proof of part (ii) of Theorem \ref{P6thm1}. By a use of $(\ref{P6eq5})$, we get
$$
\phi_n(\widehat{f}(\pi_n)\boldsymbol{\cdot}\pi_n(x)\boldsymbol{\cdot}R_n)
\leq d_{\pi_n}^2\sqrt{\mbox{tr}(\widehat{f}(\pi_n)^*\widehat{f}(\pi_n))}\phi_n(R_n).
$$
Further, by applying the Cauchy-Schwarz inequality and making use of $(\ref{P6eq24})$, we have
\begin{align*}
&\sum_{n\in\Lambda} d_{\pi_n}\phi_n(\widehat{f}(\pi_n)\boldsymbol{\cdot}\pi_n(x)\boldsymbol{\cdot}R_n)\\
&\leq\sum_{n\in\Lambda}\sqrt{d_{\pi_n}\mbox{tr}(\widehat{f}(\pi_n)^*\widehat{f}(\pi_n))}\sqrt{d_{\pi_n}^{5}}\phi_n(R_n)
\\
&\leq\sqrt{\sum_{n\in\Lambda} d_{\pi_n}\mbox{tr}(\widehat{f}(\pi_n)^*\widehat{f}(\pi_n))}\sqrt{\sum_{n\in\Lambda} d_{\pi_n}^5\phi_n^2(R_n)}\\
&\leq\sqrt{1-\widehat{f}(\pi_0)^2}\sqrt{\sum_{n\in\Lambda} d_{\pi_n}^5\phi_n^2(R_n)}\,,
\end{align*}
and the final step of the proof follows from a direct computation using the above inequality.

Now to prove part (iii), same ingredients as in the proof of part (iii) of Theorem \ref{P6thm1} will be used, i.e. Theorem \ref{TheoE}, Lemma \ref{P6lem1}, $\sigma_i(\pi_n(x))=1$, the Cauchy-Schwarz inequality, $(\ref{P6eq10})$, $(\ref{P6eq11})$ and finally Theorem \ref{TheoC}(b) are to be used, one by one.
That yields the following:
\begin{align*}
&|M_{\chi_n}(\widehat{f}(\pi_n)\boldsymbol{\cdot}\pi_n(x)\boldsymbol{\cdot}R_n)|^2\\
&\leq\frac{1}{d_{\pi_n}}\sum_{i=1}^{d_{\pi_n}}\sigma_i^{2d_{\pi_n}}(\widehat{f}(\pi_n)\boldsymbol{\cdot}\pi_n(x)\boldsymbol{\cdot}R_n)\\
&\leq\frac{1}{d_{\pi_n}}\left(\sum_{i=1}^{d_{\pi_n}}\sigma_i^{d_{\pi_n}}(\widehat{f}(\pi_n)\boldsymbol{\cdot}\pi_n(x)\boldsymbol{\cdot}R_n)\right)^2\\
&\leq\frac{1}{d_{\pi_n}}\left(\sum_{i=1}^{d_{\pi_n}}\sigma_i^{d_{\pi_n}}(\widehat{f}(\pi_n))\sigma_i^{d_{\pi_n}}(\pi_n(x))\sigma_i^{d_{\pi_n}}(R_n)\right)^2\\
&=\frac{1}{d_{\pi_n}}\left(\sum_{i=1}^{d_{\pi_n}}\sigma_i^{d_{\pi_n}}(\widehat{f}(\pi_n))\sigma_i^{d_{\pi_n}}(R_n)\right)^2\\
&\leq\frac{1}{d_{\pi_n}}\left(\sum_{i=1}^{d_{\pi_n}}\sigma_i^{2d_{\pi_n}}(\widehat{f}(\pi_n))\right)
\left(\sum_{i=1}^{d_{\pi_n}}\sigma_i^{2d_{\pi_n}}(R_n)\right)\\
&\leq\left(\sum_{i=1}^{d_{\pi_n}}\sigma_i^2(\widehat{f}(\pi_n))\right)^{d_{\pi_n}}\hspace{-10pt}\sigma_1^{2d_{\pi_n}}(R_n)\\
&\leq\left(\mbox{tr}(\widehat{f}(\pi_n)^*\widehat{f}(\pi_n))\right)^{d_{\pi_n}}\hspace{-5pt}\left(\frac{\|R_n\|_{\pi_n}}{\|E_{11}^{(n)}\|_{\pi_n}}\right)^{2d_{\pi_n}}.
\end{align*}
In other words,
$$
|M_{\chi_n}(\widehat{f}(\pi_n)\boldsymbol{\cdot}\pi_n(x)\boldsymbol{\cdot}R_n)|^{\frac{1}{d_{\pi_n}}}
\leq\sqrt{\mbox{tr}(\widehat{f}(\pi_n)^*\widehat{f}(\pi_n))}\left(\frac{\|R_n\|_{\pi_n}}{\|E_{11}^{(n)}\|_{\pi_n}}\right).
$$
Then, using the Cauchy-Schwarz inequality and $(\ref{P6eq24})$, we get
\begin{align*}
&\sum_{n\in\Lambda} d_{\pi_n}|M_{\chi_n}(\widehat{f}(\pi_n)\boldsymbol{\cdot}\pi_n(x)\boldsymbol{\cdot}R_n)|^{\frac{1}{d_{\pi_n}}}\\
&\leq\sum_{n\in\Lambda}\sqrt{d_{\pi_n}\mbox{tr}(\widehat{f}(\pi_n)^*\widehat{f}(\pi_n))}\left(\sqrt{d_{\pi_n}}\frac{\|R_n\|_{\pi_n}}{\|E_{11}^{(n)}\|_{\pi_n}}\right)\\
&\leq\sqrt{\sum_{n\in\Lambda} d_{\pi_n}\mbox{tr}(\widehat{f}(\pi_n)^*\widehat{f}(\pi_n))}
\sqrt{\sum_{n\in\Lambda} d_{\pi_n}\frac{\|R_n\|_{\pi_n}^2}{\|E_{11}^{(n)}\|_{\pi_n}^2}}\\
&\leq\sqrt{1-\widehat{f}(\pi_0)^2}\sqrt{\sum_{n\in\Lambda}  d_{\pi_n}\frac{\|R_n\|_{\pi_n}^2}{\|E_{11}^{(n)}\|_{\pi_n}^2}}.
\end{align*}
Rest of the proof relies on a straightforward computation.

Now if $\widehat{f}(\pi_0)=1$, from $(\ref{P6eq24})$ we have $\mbox{tr}(\widehat{f}(\pi_n)^*\widehat{f}(\pi_n))=0=\sum_{i=1}^{d_{\pi_n}}\sigma_i^2(\widehat{f}(\pi_n))$ for any $n\in\Lambda$, i.e. $\sigma_i(\widehat{f}(\pi_n))=0$ for $1\leq i\leq d_{\pi_n}$, or equivalently $\widehat{f}(\pi_n)=0$, $0$ being the zero matrix in $\M_{d_{\pi_n}}(\IC)$. Clearly, in this case equality holds at $(\ref{P6eq32})$, $(\ref{P6eq33})$ and $(\ref{P6eq34})$ for any arbitrary choice of \lq$R_n$\rq s. The proof of Theorem \ref{P6thm2} is therefore complete.
\epf

At this point, we take $G$ to be abelian. As a consequence,  $d_{\pi_n}=1$ for all $n\in\Lambda$ (cf. \cite[Corollary (3.6), p. 71]{Foll}), and hence $\IH_{\pi_n}$ is taken to be $\IC$ for all $n\in\Lambda$. Therefore the unitary $\pi_n(x)\in\M_1(\IC)$ for each $x\in G$, i.e. $\pi_n$ is precisely a continuous homomorphism from $G$ into $\mathbb{T}$.
Now for any $f\in L^\infty(G, \mathcal{B}(\IH))\subset L^1(G, \mathcal{B}(\IH))$ ($\IH$ is any complex Hilbert space), the Fourier transform $\widehat{f}(\pi_n)$ can be defined in exactly the same way like $(\ref{P6eq6})$ (the integral is understood in the sense of Bochner). We further define a subspace $F^\infty(G, \mathcal{B}(\IH))$ of $L^\infty(G, \mathcal{B}(\IH))$ by
$$
F^\infty(G, \mathcal{B}(\IH))
=\{f\in L^\infty(G, \mathcal{B}(\IH)):\widehat{f^*}(\pi_n)=0, n\in\Lambda\}.
$$
Here $f^*(x):=(f(x))^*$ for each $x\in G$. The goal of the following theorem is to characterize the $p$-uniform $\IC$-convexity of the space $\mathcal{B}(\IH)$ in terms of the Bohr phenomenon for the members of $F^\infty(G, \mathcal{B}(\IH))$.
\bthm\label{P6thm3}
The Banach space $\IB(\IH)$ with the operator norm $\|.\|_{\IB(\IH)}$ is $p$-uniformly $\IC$-convex \emph{($2\leq p<\infty$)} if and only if
there exists a constant $r_0(\IB(\IH))>0$ such that
\be\label{P6eq37}
\left(\|\widehat{f}(\pi_0)\|_{\IB(\IH)}^p+\sum_{n\in\Lambda} r^n\|\widehat{f}(\pi_n)\|_{\IB(\IH)}^p\right)^{\frac{1}{p}}
\leq\|f\|_{L^\infty(G, \mathcal{B}(\IH))}
\ee
for all $0\leq r\leq r_0(\IB(\IH))$ and for all $f\in F^\infty(G, \mathcal{B}(\IH))$, where $G$ is a compact abelian second countable group.
\ethm
\bpf
We prove the ``if" part first.
Since $|\pi_n(x)|=1$ for all $x\in G$,
for any $n\in\Lambda$
there exists a function $\theta_n:G\to\IR$ such that $\pi_n(x)=e^{i\theta_n(x)}, x\in G$. Now given any $A, B\in\IB(\IH)$, we consider $f(x)=A+\pi_{1}(x)B\in F^\infty(G, \IB(\IH))$.
An appropriate use of the Peter-Weyl theorem gives that $\widehat{f}(\pi_0)=A$, $\widehat{f}(\pi_1)=B$ and $\widehat{f}(\pi_n)=0$ for all $n\in\Lambda\setminus\{1\}$. According to our hypothesis, from $(\ref{P6eq37})$ we have
$$
\left(\|A\|_{\IB(\IH)}^p+r_0(\IB(\IH))\|B\|_{\IB(\IH)}^p\right)^{\frac{1}{p}}
\leq\sup_{x\in G}\|A+e^{i\theta_1(x)}B\|_{\IB(\IH)}
=\sup_{y\in S}\|A+yB\|_{\IB(\IH)},
$$
where $S=\{y=e^{i\theta_1(x)}:x\in G\}\subset\mathbb{T}=\{e^{i\theta}:\theta\in\IR\}$,
and hence
$$
\sup_{y\in S}\|A+yB\|_{\mathcal{B}(\IH)}
\leq\sup_{y\in \mathbb{T}}\|A+yB\|_{\mathcal{B}(\IH)}
=\max_\theta\|A+e^{i\theta}B\|_{\IB(\IH)}.
$$
Combining the above two inequalities, we conclude that
$\IB(\IH)$ is $p$-uniformly $\IC$-convex.

We prove the converse now. Given any $f\in F^\infty(G,\IB(\IH))$, without loss of generality we assume that $\|f\|_{L^\infty(G,\IB(\IH))}=1$. Then we construct $f_{u,v}:G\to\IC$, $u, v\in\IH$ with $\|u\|_\IH=\|v\|_\IH=1$, given by
$$
f_{u,v}(x)=\langle f(x)u, v\rangle
\in L^1(G, \IC).
$$
We clarify that $\langle .,.\rangle$ means the inner product on $\IH$ in this theorem too, and $\|.\|_{\IH}$ is the norm on $\IH$ induced by this inner product.
We further choose $\theta(u, v)$ such that $e^{i\theta(u, v)}\widehat{f_{u,v}}(\pi_0)\geq 0$, and write $F_{u,v}(x)=e^{i\theta(u, v)}f_{u,v}(x)$.
Thus
$$|F_{u,v}(x)|=|\langle f(x)u, v\rangle|\leq\|f(x)\|_{\IB(\IH)}\leq 1 \,\mbox{a.e. for}\, x\in G,
$$
thereby implying $F_{u,v}\in L^\infty(G,\IC)\subset L^1(G,\IC)$ and $\mbox{Re}(F_{u,v})\leq 1$ a.e. on $G$. Also, $\widehat{F_{u,v}}(\pi_0)\geq 0$. From $(\ref{P6eq2})$, it follows that
$$
\left|\widehat{F_{u,v}}(\pi_n)+\ov{\widehat{F_{u,v}}(\ov{\pi_n})}
\right|
\leq 2\left(1-\widehat{F_{u,v}}(\pi_0)\right)
$$
for any $n\in\Lambda$.
As $\ov{\widehat{F_{u,v}}(\ov{\pi_n})}=e^{-i\theta(u,v)}\langle\widehat{f^*}(\pi_n)v, u\rangle=0$, the above inequality yields
$$
|\widehat{F_{u,v}}(\pi_n)|\leq 2(1-\widehat{F_{u,v}}(\pi_0)).
$$
As a consequence,
\begin{align*}
&\|\widehat{f}(\pi_0)+(e^{i\theta}/2)\widehat{f}(\pi_n)\|_{\IB(\IH)}\\
&=\sup_{\substack{\|u\|_{\IH}=1\\\|v\|_{\IH}=1}}|\langle(\widehat{f}(\pi_0)+(e^{i\theta}/2)\widehat{f}(\pi_n))u, v\rangle|\\
&=\sup_{\substack{\|u\|_{\IH}=1\\\|v\|_{\IH}=1}}|\widehat{F_{u,v}}(\pi_0)+(e^{i\theta}/2)\widehat{F_{u,v}}(\pi_n)|\\
&\leq\sup_{\substack{\|u\|_{\IH}=1\\\|v\|_{\IH}=1}}\widehat{F_{u,v}}(\pi_0)+\frac{1}{2}\sup_{\substack{\|u\|_{\IH}=1\\\|v\|_{\IH}=1}}|\widehat{F_{u,v}}(\pi_n)|\\
&\leq 1
\end{align*}
for any $\theta\in\IR$. Since $\IB(\IH)$ is $p$-uniformly $\IC$-convex, there exists $\lambda(\IB(\IH))>0$ such that
$$
\|\widehat{f}(\pi_0)\|_{\IB(\IH)}^p+(\lambda(\IB(\IH))/2)\|\widehat{f}(\pi_n)\|_{\IB(\IH)}^p\leq
\left(\max_\theta\|\widehat{f}(\pi_0)+(e^{i\theta}/2)\widehat{f}(\pi_n)\|_{\IB(\IH)}\right)^p\leq 1.
$$
Hence
$$
\|\widehat{f}(\pi_n)\|_{\IB(\IH)}^p
\leq \frac{2}{\lambda(\IB(\IH))}
\left(1-\|\widehat{f}(\pi_0)\|_{\IB(\IH)}^p
\right)
$$
for any $n\in\Lambda$. From here, it is evident that $(\ref{P6eq37})$ holds for any nonnegative real number $r$ satisfying
$$
r\leq\lambda(\IB(\IH))/(2+\lambda(\IB(\IH))).
$$
Our proof is therefore complete.
\epf

\section{Some related remarks}
We end this article with a few relevant comments and observations.
\bee
\item
If we set $G=\mathbb{T}$, i.e. the multiplicative group of complex numbers with unit modulus, then $d_{\pi_n}=1$ for all $n\in\Lambda$, and $\mathcal{H}_{\pi_n}$ can be taken to be $\IC$ for all $n\in\Lambda$. It is known that
$\widehat{\mathbb{T}}\cong\mathbb{Z}$ (cf. \cite[Theorem (4.5), p. 89]{Foll}), the irreducible unitary representation corresponding to $n\in\mathbb{Z}$ is given by $x^n$, $x\in\mathbb{T}$,
and the Fourier transform of any $f\in L^1(\mathbb{T}, \IC)$ associated with this representation is denoted by $\widehat{f}(n)$.
We therefore take $\Lambda=\IN$,
$\pi_0(x)=1$ (as usual),
$\pi_n(x)=x^{k}$ for all $n=2k-1, k\geq 1$ and $\pi_n(x)=x^{-k}$ for all $n=2k, k\geq1$, where $x=e^{i\theta}\in\mathbb{T}$. Also take $R_n=z^{k}$ for all $n=2k-1$ and $R_n=0$ for all $n=2k$,
$z\in\IC$ and $k\geq 1$.
Further, for all $n\in\IN$, $\|.\|_{\pi_n}$, $\phi_n$ and $M_{\chi_n}$ are taken to be the usual absolute value of complex numbers, and the $\boldsymbol{\cdot}$ product is taken to be the usual product of complex numbers. Finally, we assume that $f\in L^1(\mathbb{T}, \IC)$ has analytic Fourier series, i.e. $\widehat{f}(\pi_{2k})=0$ for all $k\in\mathbb{N}$. Using the fact that $\ov{\pi_{2k-1}}(x)=\pi_{2k}(x)$, we observe that under all the above assumptions, all three inequalities $(\ref{P6eq29})$, $(\ref{P6eq30})$ and $(\ref{P6eq31})$ from Theorem \ref{P6thm1} convert to
\begin{align*}
\hspace{20pt}&\widehat{f}(\pi_0)+\sum_{n=1}^\infty\left|\left(\widehat{f}(\pi_{n})+\ov{\widehat{f}(\ov{\pi_n})}\right)\pi_{n}(x)R_{n}
\right|\\
&=\widehat{f}(\pi_0)+\sum_{k=1}^\infty|\widehat{f}(\pi_{2k-1})\pi_{2k-1}(x)R_{2k-1}|\\
&=\widehat{f}(0)+\sum_{k=1}^\infty|\widehat{f}(k)e^{ik\theta}|r^k
\leq 1,
\end{align*}
where $|z|=r$. According to the parts (i), (ii) and (iii) of Theorem \ref{P6thm1}, if $f$ satisfies all the hypotheses of this theorem, then the above inequality is satisfied whenever
\begin{align*}
&\sum_{n=1}^\infty\|R_n\|_{\pi_n}
=\sum_{n=1}^\infty d_{\pi_n}^2\phi_n(R_n)
=\sum_{n=1}^\infty\left(\|R_n\|_{\pi_n}/\|E_{11}^{(n)}\|_{\pi_n}\right)\\
&=\sum_{k=1}^\infty|R_{2k-1}|
=\sum_{k=1}^\infty|z^k|\\
&=r/(1-r)
\leq1/2,
\end{align*}
i.e. whenever $r\leq 1/3$. Thus we get Bohr's original theorem, i.e. Theorem \ref{TheoA} back here. Similarly, if $f$ satisfies all the hypotheses of Theorem \ref{P6thm2}, then from each of the three parts of this theorem we get that the inequality
$$
\widehat{f}(0)+\sum_{k=1}^\infty|\widehat{f}(k)e^{ik\theta}|r^k
\leq 1
$$
holds if
$$
\left(1+\widehat{f}(0)\right)\left(1+\sum_{k=1}^\infty r^{2k}\right)
\leq 2,
$$
i.e. if
$r\leq\sqrt{(1-\widehat{f}(0))/2}$, which reduces to $r\leq 1/\sqrt{2}$ for $\widehat{f}(0)=0$. In other words, a part of \cite[Corollary 2.9]{Paul1} that improves on Theorem \ref{TheoA} (for analytic self mappings of $\ID$ fixing the origin) has been recovered here.
\item
It should be noted that in the articles \cite{Abu, Pon}, Bohr inequalities analogous to the inequalities in Theorem \ref{P6thm1} and Theorem \ref{P6thm2} were established for complex valued harmonic functions defined in $\ID$.
\item
Looking at Theorem \ref{P6thm1}, one can ask if it is possible to exclude the term $\ov{\widehat{f}(\ov{\pi_n})}$ from the Bohr inequalities. Modifying an example from \cite{Abu}, we show that the answer is negative in general. Namely, we define $f_\mu:\mathbb{T}\to\IC$ by
$$
f_\mu(e^{i\theta})= \cos \mu+ i \sin \mu\left( \sum_{n\in\mathbb{Z}\setminus\{0\}}\frac{1}{4n^2}e^{in\theta}\right)
=\sum_{n=-\infty}^\infty a_n(\mu)e^{in\theta}
$$
where $\theta\in[0, 2\pi)$ and $\mu\in(0,\pi/2)$. Clearly $|f_\mu(e^{i\theta})|<1$ and hence $\mbox{Re}(f_\mu)<1$ on all of $\mathbb{T}$. Also, $a_0(\mu)> 0$.
Now suppose that for all $\mu$, there exists a constant $r_0>0$ such that
$$
a_0(\mu)+\sum_{n\in\mathbb{Z}\setminus\{0\}}|a_n(\mu)e^{in\theta}z_n|\leq 1
$$
\,holds for any given sequence $\{z_n\}_{n\in\mathbb{Z}\setminus\{0\}}$ of complex numbers satisfying
$\sum_{n\in\mathbb{Z}\setminus\{0\}}|z_n|\leq r_0$, and $z_{\tilde{n}}\neq 0$ for at least one $\tilde{n}\in\mathbb{Z}\setminus\{0\}$. Thus we have
$$
\frac{1}{4\tilde{n}^2}\left(\frac{\sin\mu}{1-\cos\mu}\right)\leq \frac{1}{|z_{\tilde{n}}|}
$$
for any $f_\mu$ as above. Letting $\mu\to 0$, we see that the left hand side goes to infinity while the right hand side is a constant, which leads to a contradiction.
\eee


\begin{thebibliography}{99}
\bibitem{Abu} {\sc Y. Abu Muhanna}: Bohr's phenomenon in subordination and bounded harmonic classes,
\textit{Complex Var. Elliptic Equ.}, 55 (2010), no. 11, 1071--1078.
\bibitem{Aiz} {\sc L. Aizenberg}: Multidimensional analogues of Bohr's theorem on power series, \textit{Proc. Amer. Math. Soc.}, 128 (2000), no. 4, 1147--1155.
\bibitem{Ayt} {\sc A. Aytuna, P. Djakov}: Bohr property of bases in the space of entire functions and its generalizations,
\textit{Bull. Lond. Math. Soc.}, 45 (2013), no. 2, 411--420.
\bibitem{Bala} {\sc R. Balasubramanian, B. Calado, H. Queff{\'e}lec}: The Bohr inequality for ordinary Dirichlet series,
\textit{Studia Math.}, 175 (2006), no. 3, 285--304.
\bibitem{Bay} {\sc F. Bayart, D. Pellegrino, J. B. Seoane-Sep{\'u}lveda}: The Bohr radius of the $n$-dimensional
polydisk is equivalent to $\sqrt{(\log n)/n}$, \textit{Adv. Math.}, 264 (2014), 726--746.
\bibitem{Bla} {\sc O. Blasco, M. Pavlovi{\'c}}: Complex convexity and vector-valued Littlewood-Paley inequalities, \textit{Bull. London Math. Soc.}, 35 (2003), no. 6, 749--758.
\bibitem{Bla2} {\sc O. Blasco}: The $p$-Bohr radius of a Banach space, \textit{Collect. Math.}, 68 (2017), no. 1, 87--100.
\bibitem{Boas} {\sc H. P. Boas, D. Khavinson}: Bohr's power series theorem in several variables,
\textit{Proc. Amer. Math. Soc.}, 125 (1997), no. 10, 2975--2979.
\bibitem{Bohr} {\sc H. Bohr}: A theorem concerning power series,
\textit{Proc. London Math. Soc.}, (2) 13 (1914), 1--5.
\bibitem{Def2} {\sc A. Defant, D. Garc{\'i}a, M. Maestre, D. P{\'e}rez-Garc{\'i}a}: Bohr's strip for vector valued Dirichlet series,
\textit{Math. Ann.}, 342 (2008), no. 3, 533--555.
\bibitem{Dix} {\sc P. G. Dixon}: Banach algebras satisfying the non-unital von Neumann inequality,
\textit{Bull. London Math. Soc.}, 27 (1995), no. 4, 359--362.
\bibitem{Foll} {\sc G. B. Folland}: A course in abstract harmonic analysis. Studies in Advanced Mathematics, \textit{CRC Press, Boca Raton, FL}, 1995.
\bibitem{Ga} {\sc D. Galicer, M. Mansilla, S. Muro}: Mixed Bohr radius in several variables, \textit{Trans. Amer. Math. Soc.}, 373 (2020), no. 2, 777--796.
\bibitem{Cue} {\sc J. Garcia-Cuerva, K. S. Kazarian, V. I. Kolyada, and J. L. Torrea}: Vector-valued Hausdorff-Young inequality and applications, \textit{Russian Math. Surveys}, 53 (1998), no. 3, 435--513.
\bibitem{Ham} {\sc H. Hamada, T. Honda, G. Kohr}: Bohr's theorem for holomorphic mappings with values in homogeneous balls,
\textit{Israel J. Math.}, 173 (2009), 177--187.
\bibitem{Hew} {\sc E. Hewitt, K. A. Ross}: Abstract harmonic analysis. Vol. II: Structure and analysis for compact groups. Analysis on locally compact Abelian groups. Die Grundlehren der mathematischen Wissenschaften, Band 152 \textit{Springer-Verlag, New York-Berlin}, 1970.
\bibitem{Horn1} {\sc R. A. Horn, C. R. Johnson}: Matrix analysis. Corrected reprint of the 1985 original, \textit{Cambridge University Press, Cambridge}, 1990.
\bibitem{Horn} {\sc R. A. Horn, C. R. Johnson}: Topics in matrix analysis, \textit{Cambridge University Press, Cambridge}, 1991.
\bibitem{Pon} {\sc I. R. Kayumov, S. Ponnusamy}: Bohr's inequalities for the analytic functions with lacunary series and harmonic functions, \textit{J. Math. Anal. Appl.}, 465 (2018), no. 2, 857--871.
\bibitem{La} {\sc P. Lass\`{e}re, E. Mazzilli}: Estimates for the Bohr radius of a Faber-Green condenser in the complex plane, \textit{Constr. Approx.}, 45 (2017), no. 3, 409--426.
\bibitem{Li} {\sc C.-K. Li, R. Mathias}: Generalizations of Ky Fan's dominance Theorem, \textit{SIAM J. Matrix Anal. Appl.}, 19 (1998), no. 1, 99--106.
\bibitem{Marc} {\sc M. Marcus, H. Minc}: Generalized matrix functions, \textit{Trans. Amer. Math. Soc.}, 116 (1965), 316--329.
\bibitem{Paul1} {\sc V. I. Paulsen, G. Popescu, D. Singh}: On Bohr's inequality,
\textit{Proc. London Math. Soc.}, (3) 85 (2002), no. 2, 493--512.
\bibitem{Paul2} {\sc V. I. Paulsen, D. Singh}: Bohr's inequality for uniform algebras, \textit{Proc. Amer. Math. Soc.}, 132 (2004), no. 12, 3577--3579.
\bibitem{Paul3} {\sc V. I. Paulsen, D. Singh}: Extensions of Bohr's inequality, \textit{Bull. London Math. Soc.}, 38 (2006), no. 6, 991--999.
\bibitem{Pop1} {\sc G. Popescu}: Bohr inequalities for free holomorphic functions on polyballs,
\textit{Adv. Math.}, 347 (2019), 1002--1053.

\end{thebibliography}
\end{document}